\documentclass[12pt]{article}
\usepackage{amsmath,amsfonts}
\usepackage{amssymb}
\usepackage{amsthm}
\usepackage{times,txfonts}

\usepackage{accents, kbordermatrix}

\newtheorem{thm}{Theorem}[section]

\newtheorem{lemma}[thm]{Lemma}

\newtheorem{remark}[thm]{Remark}
\newtheorem{example}[thm]{Example}

\numberwithin{equation}{section}

\def\pf{\noindent{\it Proof.} \ }

\def\qed{\hfill $\square$}

\title{On a fundamental system of solutions of a certain hypergeometric equation}

\author{Teruhisa Tsuda \\
Department of Economics, Hitotsubashi University,  
\\
Tokyo 186-8601, Japan.
\\
tudateru@econ.hit-u.ac.jp}

\date{March 14, 2014\\
Revised: August 2, 2014}

\begin{document}
\maketitle

\renewcommand{\thefootnote}{\fnsymbol{footnote}}
\footnotetext{{\it 2010 Mathematics Subject Classification} 
33C65,   
33C70, 
34M56. 
} 
\footnotetext{{\it Keywords:} 
hypergeometric function, 
isomonodromic deformation, 
Painlev\'e equation.}

\begin{abstract}
We study the linear Pfaffian systems 
satisfied by 
a certain class of hypergeometric functions, which includes
Gau\ss's ${}_2 F_{1}$,
Thomae's ${}_L F_{L-1}$  
and 
Appell--Lauricella's $F_D$.
In particular, we 
present
a fundamental system of solutions 
with a characteristic local behavior
by means of 
Euler-type integral representations.
We also 
discuss 
how they are related to the theory of isomonodromic deformations
or Painlev\'e equations.
\end{abstract}

\section{Introduction}

\subsection{Hypergeometric function $F_{L,N}$}

Fix integers $L \geq 2$ and $N \geq 1$.
We consider the hypergeometric function
 $F_{L,N}=F_{L,N}(\alpha, \beta, \gamma; x)$
 in $N$ variables $x=(x_1, \ldots, x_N)$
defined by means of the power series 
\begin{equation}
\label{eq:hgseries}
F_{L,N}(\alpha, \beta, \gamma; x)
=
\sum_{m_i \geq0}
\frac{(\alpha_1)_{|m|} \cdots (\alpha_{L-1})_{|m|} (\beta_1)_{m_1} \cdots  (\beta_N)_{m_N} }{ (\gamma_1)_{| m|} \cdots (\gamma_{L-1})_{| m|} (1)_{m_1} \cdots (1)_{m_N}} 
{x_1}^{m_1} \cdots {x_N}^{m_N}
\end{equation}
convergent in the polydisc 
$D_0=\{|x_1|<1,\ldots,|x_N|<1\} \subset {\mathbb C}^N$.
Here 
$|m|=m_1+ \cdots +m_N$
and
$(a)_n=\Gamma(a+n)/\Gamma(a)$
is the Pochhammer symbol.
The $2L+N-2$ parameters
\[
(\alpha,\beta,\gamma)=(\alpha_1,\ldots,\alpha_{L-1}, \beta_1,\ldots,\beta_N, \gamma_1, \ldots,\gamma_{L-1})
\]
are complex constants such that 
$\gamma_n \notin {\mathbb Z}_{< 0}$.
Note that,
if $(L,N)=(2,1)$, $(L,1)$
and
$(2,N)$, then
the hypergeometric function $F_{L,N}$ reduces  
to Gau\ss's ${}_2F_1$, 
Thomae's ${}_L F_{L-1}$ \cite{tho70} 
and 
Appell--Lauricella's $F_D$ \cite{ak26, lau93},
respectively.

It is straightforward to verify from the power series
 (\ref{eq:hgseries})
that $F_{L,N}$ solves the system 
of linear differential equations
\begin{equation} \label{eq:hgeq}
\left\{
x_i  \left( \beta_i+\delta_i\right)
\prod_{k=1}^{L-1} \left( \alpha_k+ 
{\cal D}
\right) 
- \delta_i
 \prod_{k=1}^{L-1} \left( \gamma_k-1+  {\cal D} \right) 
\right\}y=0
\quad
(1 \leq i \leq N)
\end{equation}
where 
\[
\delta_i=x_i \frac{\partial}{\partial x_i}
\quad \text{and} \quad
{\cal D}=\sum_{i=1}^N \delta_i.
\]
Moreover, $F_{L,N}$ possesses 
an Euler-type integral representation (see \cite[Proposition~2.1]{tsu12}) 
\[
F_{L,N}
= 
\prod_{k=1}^{L-1} \frac{\Gamma(\gamma_k)}{\Gamma(\alpha_k) \Gamma(\gamma_k-\alpha_k)}
\times
\int_{\Delta_0}
U(t) \varphi_0
\]
with the domain $\Delta_0$ being an 
$(L-1)$-simplex
\begin{equation} \label{eq:Delta_0}
\Delta_0=\{0 \leq t_{L-1} \leq \cdots \leq t_2 \leq t_1 \leq 1 \} \subset {\mathbb R}^{L-1}
\end{equation}
and the integrand $U(t) \varphi_0$ given by (\ref{eq:U}) and (\ref{eq:varphi}) below.
Starting from the above integral representation,
we have a certain linear Pfaffian system
 (see \cite[Theorem~2.2]{tsu12})
equivalent to (\ref{eq:hgeq}), 
which is the main object studied in this paper.
We shall briefly
introduce this linear Pfaffian system, denoted by ${\cal P}_{L,N}$, 
in the following Sect.~\ref{subsect:pfaff}.

\subsection{Linear Pfaffian system ${\cal P}_{L,N}$ (the hypergeometric equation)}
\label{subsect:pfaff}

Let 
\[
\zeta_k=
\alpha_k- \gamma_{k+1}  
\quad (\gamma_L=1), \quad 
\eta_k=\gamma_k-\alpha_k,
\quad
\theta_i=-\beta_i
\]
and 
consider a multi-valued function 
\begin{equation} \label{eq:U}
U(t)=
\prod_{k=1}^{L-1} {t_k}^{\zeta_k}  
(t_{k-1}-t_k)^{\eta_k}
\prod_{i=1}^{N} (1-x_i t_{L-1})^{\theta_i}
\end{equation}
in $t=(t_1,t_2,\ldots,t_{L-1})$
with $t_0=1$.
Consider the rational $(L-1)$-forms
\begin{align}
\varphi_0&= \frac{\underline{{\rm d}t}}{ \prod_{k=1}^{L-1}(t_{k-1}-t_k)}, 
\nonumber
\\
\varphi_n^{(i)}&= \frac{\underline{{\rm d}t}}{ (x_i t_{L-1}-1)\prod_{\begin{subarray}{l} k=1 \\ k \neq n \end{subarray}}^{L-1}(t_{k-1}-t_k)}
\quad 
\left(
\begin{array}{c}
 1 \leq i \leq N \\
1 \leq n \leq L-1
\end{array}
\right)
\label{eq:varphi}
\end{align}
where
\begin{equation} \label{eq:dt}
\underline{{\rm d}t}
= {\rm d}t_1 \wedge \cdots \wedge {\rm d} t_{L-1}.
\end{equation}
Now we define the vector-valued function
\[\boldsymbol{y}=
\boldsymbol{y}(x; \Delta)=
{}^{\rm T}\left(
y_0,y_1^{(1)},\ldots,y_{L-1}^{(1)},
y_1^{(2)},\ldots,y_{L-1}^{(2)},
 \ldots,
y_1^{(N)},\ldots,y_{L-1}^{(N)}
\right)
\]
by the integrals
\begin{equation} \label{eq:int}
y_0=\int_{\Delta} U(t) \varphi_0,
\quad
y_n^{(i)}=\int_{\Delta} U(t) \varphi_n^{(i)}
\end{equation}
over a suitable domain $\Delta$.
The function
$\boldsymbol{y}$ then satisfies the linear Pfaffian system
\[
\tag{${\cal P}_{L,N}$}
{\rm d} \boldsymbol{y}
= \left\{
\sum_{i=1}^N \left(E_i {\rm d} \log x_i+F_i  {\rm d} \log (1-x_i)
\right)+
\sum_{1 \leq i< j \leq N} G_{ij}   {\rm d} \log (x_i-x_j)
\right\}
\boldsymbol{y}
\]
of rank $N(L-1)+1$.
Here, the coefficient matrices are linear functions in the constant parameters
$(\alpha,\beta,\gamma)$
given by
\begin{allowdisplaybreaks}
\begin{align*}
&E_i=
\kbordermatrix{
&0 \hspace{-2mm}&&1&& 2 &&    &&&&i&&&&       &&\hspace{-2mm} N
\\
&& \vrule && \vrule && \vrule && \vrule &&&&&&\vrule&& \vrule
\\
\cline{2-18}
&  a_1 \hspace{-2mm}&\vrule&  &\vrule&  &\vrule&  &\vrule& \hspace{-2mm}  b_{i,1}&&&&&\vrule&&\vrule& 
\\
& a_2 \hspace{-2mm} &\vrule&  &\vrule&  &\vrule&  &\vrule& \hspace{-2mm}  a_2&b_{i,2} &&&&\vrule&&\vrule& 
\\
 i & a_3 \hspace{-2mm} &\vrule& \hspace{-2mm} -\beta_1 I_{L-1} \hspace{-2mm} &\vrule& \hspace{-2mm} -\beta_2 I_{L-1} \hspace{-2mm}&\vrule& \hspace{-2mm}\cdots \hspace{-2mm} &\vrule& \hspace{-2mm} a_3&a_3& b_{i,3} &&&\vrule&\hspace{-2mm} \cdots \hspace{-2mm} &\vrule&\hspace{-2mm}  -\beta_N I_{L-1} 
 \\
 &\vdots \hspace{-2mm} &\vrule&  &\vrule&  &\vrule&  &\vrule&\hspace{-2mm}  \vdots&\vdots&\ddots&\ddots&&\vrule&&\vrule&  
 \\
 & a_{L-1} \hspace{-2mm} &\vrule&  &\vrule&  &\vrule&  &\vrule& \hspace{-2mm} a_{L-1}& a_{L-1} &\cdots&a_{L-1}&b_{i,L-1} \hspace{-2mm}  &\vrule&&\vrule&  
\\
\cline{2-18}
&& \vrule && \vrule && \vrule && \vrule &&&&&&\vrule&& \vrule
},
\\
&F_i=
\kbordermatrix{
&0\hspace{-2mm}&& && & i & && \\
0&-\beta_i \hspace{-2mm} &\vrule& &\vrule& \hspace{-2mm}  -\beta_i & -\beta_i & \cdots& -\beta_i \hspace{-2mm}  &\vrule& \\
\cline{2-11}
&&\vrule&&\vrule&&&&&\vrule& \\
\cline{2-11}
&-a_1 \hspace{-2mm} &\vrule& &\vrule& \hspace{-2mm} -a_1 &-a_1& \cdots& -a_1 \hspace{-2mm} &\vrule& \\
i &-a_2 \hspace{-2mm} &\vrule& &\vrule&\hspace{-2mm}  -a_2& -a_2&\cdots& -a_2 \hspace{-2mm}  &\vrule& \\
&\vdots\hspace{-2mm} &\vrule&&\vrule&\hspace{-2mm}  \vdots&\vdots&\ddots&\vdots \hspace{-2mm} &\vrule& \\
& -a_{L-1} \hspace{-2mm} &\vrule& &\vrule&\hspace{-2mm}  -a_{L-1}& -a_{L-1}& \cdots&-a_{L-1}  \hspace{-2mm}  &\vrule& \\
\cline{2-11}
&&\vrule&&\vrule&&&&&\vrule& \\
}, 
\quad
G_{ij}=
\kbordermatrix{
&&& i &&&&j&& \\
&&\vrule&&\vrule&&\vrule&&\vrule&  \\
\cline{2-10}
i&& \vrule &\hspace{-2mm}  -\beta_j I_{L-1} \hspace{-2mm}  & \vrule&&\vrule&  \hspace{-2mm}  \beta_j I_{L-1} \hspace{-2mm}  &\vrule& \\
\cline{2-10}
&&\vrule&&\vrule&&\vrule&&\vrule& \\
\cline{2-10}
j&& \vrule & \hspace{-2mm}  \beta_i I_{L-1}  \hspace{-2mm}  & \vrule&&\vrule& \hspace{-2mm}  -\beta_i I_{L-1} \hspace{-2mm}  &\vrule& \\
\cline{2-10}
&&\vrule&&\vrule&&\vrule&&\vrule& \\
}
\end{align*}
\end{allowdisplaybreaks}
\\
\noindent
and $a_n=\alpha_n-\gamma_n$ and
$b_{i,n}=\sum_{j \neq i} \beta_j -\gamma_n $.
The symbol $I_{L-1}$ denotes the identity matrix of size $L-1$.
We 
wrote a square matrix $M$ of size $N(L-1)+1$
with separating it into $(N+1)^2$ blocks 
as
\[
M=\kbordermatrix{
&0 \hspace{-2mm}&&1&&     && \hspace{-2mm}N \\
0&M_{00} \hspace{-2mm} & \vrule & \hspace{-2mm}  M_{01} \hspace{-2mm} &\vrule &\hspace{-2mm} \cdots \hspace{-2mm} &\vrule&
\hspace{-2mm}  M_{0N}
\\
\cline{2-8}
1& M_{10} \hspace{-2mm} & \vrule & \hspace{-2mm}  M_{11} \hspace{-2mm} &\vrule &\hspace{-2mm} \cdots \hspace{-2mm} &\vrule&
\hspace{-2mm} 
M_{1N}
\\
\cline{2-8}
      & \vdots \hspace{-2mm} &  \vrule & \hspace{-2mm}  \vdots\hspace{-2mm}  &\vrule &\hspace{-2mm}  \ddots \hspace{-2mm}  &\vrule&
\hspace{-2mm}  \vdots
\\
\cline{2-8}
N& M_{N0} \hspace{-2mm} & \vrule &\hspace{-2mm}  M_{N1}  \hspace{-2mm} &\vrule  &\hspace{-2mm}  \cdots \hspace{-2mm}  &\vrule&
\hspace{-2mm}  M_{NN}
},
\]
where
$M_{ij}$ $(i,j \neq 0)$ is a square matrix of size $L-1$
and thus 
 $M_{00}$ is a scalar,
$M_{0j}$ $(j \neq 0)$ and $M_{i0}$ $(i \neq 0)$
are
row and  column $(L-1)$-vectors,
respectively.

The singular locus of the system
${\cal P}_{L,N}$ is a union of hyperplanes
\begin{equation}  \label{eq:Xi} 
\Xi=
\bigcup_{i=1}^N
\left(
 \{x_i=0\} \cup \{ x_i=1\} \cup \{x_i= \infty \}
 \right) 
 \cup \bigcup_{1 \leq i<j \leq N}\{ x_i=x_j\}.
\end{equation}
Therefore, the holomorphic function
$F_{L,N}$
at $x=0$ can be analytically continued along any path outside
$\Xi$.
The {\it characteristic exponents} at each divisor,  
i.e.
 the eigenvalues of each residue matrix $E_i$, $F_i$ or $G_{ij}$
are listed in the following table (Riemann scheme):
\begin{equation}
\label{eq:rsofp}
\begin{array}{|c|c|}
\hline
\text{Divisor} & \text{Characteristic exponents}  \\ \hline
x_i=0 & (b_{i,1},b_{i,2}, \ldots, b_{i,L-1}, \underbrace{0, \ldots,0}_{(N-1)(L-1)+1} )
\\
\hline
x_i=1& 
(-\beta_i-\sum_{n=1}^{L-1} a_n,  \underbrace{0, \ldots,0}_{N(L-1)} )
\\
\hline
x_i=\infty &
(\alpha_1, \alpha_2, \ldots, \alpha_{L-1},
\underbrace{\beta_i, \ldots,\beta_i}_{(N-1)(L-1)+1}  )
\\
\hline
x_i=x_j 
\ (i \neq j) &
(\underbrace{-\beta_i-\beta_j, \ldots, -\beta_i-\beta_j}_{L-1},
\underbrace{0, \ldots,0}_{(N-1)(L-1)+1}  )
\\
\hline
   \end{array}
\end{equation}

\begin{remark}\rm
\label{remark:rigid}

If we regard one specific variable $x_i$ as an independent variable and all other $x_j$'s ($j \neq i$)
as fixed constants, then ${\cal P}_{L,N}$ is a Fuchsian system of 
ordinary differential equations with respect to $x_i$.  
Notice that the sum of characteristic exponents at $x_i=0,1,\infty, x_j$ $(j \neq i)$
is equal to zero; i.e. the Fuchsian relation holds.
The {\it spectral type} of this Fuchsian system
is given by the $(N+2)$-tuple
\[
\begin{array}{ll}
\underbrace{1,1, \ldots,1}_{L-1}, (N-1)(L-1)+1
&
\text{at $x_i=0, \infty$}, 
\\
1, N(L-1)
&
\text{at $x_i=1$}, 
\\
L-1, (N-1)(L-1)+1
&
\text{at $x_i=x_j$ $(j \neq i)$} 
\end{array}
\]
of partitions of $N(L-1)+1$,
which indicates how the characteristic exponents overlap at each of the singularities.
We know, by counting Katz's index \cite{kat95}, that this Fuchsian system is {\it rigid} in the sense that its global monodromy is determined only from its local monodromy,
or from its characteristic exponents, at  each singularity.
\end{remark}

\begin{remark}\rm
\label{remark:pfaff}
The linear Pfaffian system ${\cal P}_{L,N}$
can be derived from the integral representations
(\ref{eq:int}),
with the aid of twisted de Rham theory \cite{ak94}.
For details refer to \cite{tsu12}
(in which the symbols $U(t)$, $\varphi_0$ and $\varphi_n^{(i)}$ 
are slightly different from the present ones).
\end{remark}

\subsection{Holomorphic solution at the origin}
\label{subsect:holsol}

If the domain $\Delta$ of integration is chosen to be the $(L-1)$-simplex
$\Delta_0$
(see (\ref{eq:Delta_0})),
then $\boldsymbol{y}(x;\Delta_0)$ becomes holomorphic at the origin $x=0 \in {\mathbb C}^N$.
This is the unique solution of ${\cal P}_{L,N}$ holomorphic at $x=0$
(up to multiplication by constants), and is expressible in terms of the hypergeometric function $F_{L,N}(\alpha,\beta,\gamma;x)$
as
\begin{align*}
&y_0
= c F_{L,N}, \quad
y_1^{(i)}
= \frac{\alpha_1-\gamma_1}{\gamma_1}  
c F_{L,N}(\beta_i+1,\gamma_1+1),
\\
&y_2^{(i)}
= \frac{\alpha_1(\alpha_2-\gamma_2)}{\gamma_1\gamma_2} 
c F_{L,N}(\alpha_1+1,\beta_i+1,\gamma_1+1,\gamma_2+1),
\quad \ldots
\\
&y_n^{(i)}
=\frac{\alpha_1  \cdots \alpha_{n-1}(\alpha_n-\gamma_n)}{\gamma_1 \cdots \gamma_n} c F_{L,N}(\alpha_1+1,\ldots,\alpha_{n-1}+1,\beta_i+1,\gamma_1+1,\ldots,\gamma_{n}+1),
\quad \ldots
\end{align*}
where
$c=\prod_{k=1}^{L-1}\Gamma(\alpha_k)\Gamma(\gamma_k-\alpha_k)/\Gamma(\gamma_k)$.
For notational simplicity,
we used the abbreviation 
$F_{L,N}(\beta_i+1, \gamma_1+1)$ 
to mean that among the parameters $(\alpha,\beta,\gamma)$
only the indicated ones $\beta_i$ and $\gamma_1$ are shifted by one,
and so forth.
We note that the first element $y=y_0$ of $\boldsymbol{y}$ solves (\ref{eq:hgeq}) indeed.

\subsection{Main result: fundamental system of solutions of ${\cal P}_{L,N}$}

Let $D_0$ and $D^{(j)}$ $(j=1, \ldots,N)$ be the  domains in ${\mathbb C}^N$
given by
\begin{align*}
D_0&
=\{|x_1|<1,\ldots,|x_N|<1\},
\\
D^{(j)}
&= \{ |x_j|<1  \} \cap \bigcap_{i \neq j} \{ x_i \neq x_j\}.
\end{align*}
According to the Riemann scheme (\ref{eq:rsofp}),
a fundamental system of solutions of ${\cal P}_{L,N}$ of the following form
can be expected.

\begin{thm}
\label{thm:main}
There exists a fundamental system $Y(x)$ of solutions of ${\cal P}_{L,N}$
such that
\begin{align}
\nonumber
Y(x)&=\Phi(x)
{\rm \, diag \, }
\left(
1,{x_1}^{b_{1,1}}, \ldots,{x_1}^{b_{1,L-1}}, 
{x_2}^{b_{2,1}}, \ldots,{x_2}^{b_{2,L-1}}, 
\ldots,
{x_N}^{b_{N,1}},  \ldots,{x_N}^{b_{N,L-1}}
\right),
\\
\nonumber
& \text{\footnotesize \hspace{7.5mm} $0$ 
 \hspace{9mm} $1$  \hspace{16mm} $2$ 
  \hspace{25mm} $N$}
\\
\nonumber
\Phi(x)&=
\left[
\begin{array}{c|ccc|ccc|c|ccc}
* & &&& & &&&&&
\\
\hline
*& * && && &&&&&
\\
\vdots &\vdots &\ddots& &&&& &&&
\\
* &*& \cdots & * &&& & &&&
\\
\hline
*& &&&* && &&&&
\\
\vdots &&&&\vdots &\ddots& &&&&
\\
* &&&&*& \cdots & * &&&&
\\
\hline
\vdots &&& &&&& \ddots  &&&
\\
\hline
*&  && && &&&*&&
\\
\vdots &&&&& &&&\vdots &\ddots &
\\
* &&&  &&& & &*&\cdots&*
\end{array}
\right]
\\
& \label{eq:Phi}
\qquad
+
\left[
\begin{array}{c|c|c|c|c}
R_0(x) & x_1 R^{(1)}(x) &x_2 R^{(2)}(x) & \cdots &x_N R^{(N)}(x)
\end{array}
\right],
\end{align}
where $R_0$ is a column vector and other $R^{(j)}$'s 
$(j =1 , 2, \ldots,N)$
are matrices with $L-1$ columns{\rm;}
$R_0$ is holomorphic on $D_0$ and $R_0(0)=0$ and
each $R^{(j)}$ is holomorphic on $D^{(j)}$.
\end{thm}

To be more precise, if $\Phi_0$ denotes the first term of the right-hand side  of 
{\rm(\ref{eq:Phi})}, 
each element of the $0$th column of $\Phi_0$ is a nonzero constant.
Concerning the $j$th block $(j =1 , 2, \ldots, N)$, $(*)$ represents a holomorphic function on $D^{(j)}$ whose restriction to $\{ x_j=0\}$ is not identically zero. 
In the rest of this paper, we shall prove this theorem.
Notice first that the holomorphic solution at the origin given 
in Sect.~\ref{subsect:holsol} provides the $0$th column of $Y(x)$ or of $\Phi(x)$.
Others will be explicitly constructed by a systematic use of Euler-type integral representations together with an iteration of  `cyclic' transformations;
see Theorems~\ref{thm:rep1} and  \ref{thm:rep2} for $N=1$ case
and Theorem~\ref{thm:general} for general $(L,N)$ case.

In Sect.~\ref{sect:int}, 
we prepare a cyclic expression of the hypergeometric integrals 
 (\ref{eq:int}), i.e. 
the general solution of ${\cal P}_{L,N}$.
We then construct a fundamental system of solutions, having the desired local behavior,
for $N=1$ case in Sect.~\ref{sect:N=1}
and for general $(L,N)$ case in Sect.~\ref{sect:general}.   
We also indicate 
a connection between the hypergeometric functions and isomonodromic deformations of a certain Fuchsian system in  Sect.~\ref{sect:mpd}.

\section{A cyclic expression of hypergeometric integrals}
\label{sect:int}

Let us first homogenize the integrands of the hypergeometric integrals (\ref{eq:int})
by 
introducing a set of $L$ variables $\tau=(\tau_0,\tau_1, \ldots ,\tau_{L-1})$ with
$t_n= {\tau_n}/{\tau_0}$.
Observe that
\[
\underline{{\rm d}t}={\tau_0}^{-L} \omega,
\quad \text{where} \quad
\omega=\sum_{n=0}^{L-1} (-1)^n 
\tau_n
 {\rm d}\tau_0 \wedge \cdots \wedge  \widehat{{\rm d}\tau_n}
 \wedge \cdots\wedge {\rm d} \tau_{L-1},
\]
and thereby
\begin{align*}
\varphi_0&=
\frac{\omega}{ \tau_0 \prod_{k=1}^{L-1}(\tau_{k-1}-\tau_k)},
\\
\varphi_n^{(i)}&=
\frac{\omega}{\tau_0 (x_i \tau_{L-1}-\tau_0)\prod_{\begin{subarray}{l} k=1 \\ k \neq n \end{subarray}}^{L-1}(\tau_{k-1}-\tau_k)}
\quad 
\left(
\begin{array}{c}
 1 \leq i \leq N \\
1 \leq n \leq L-1
\end{array}
\right).
\end{align*}
Cf. (\ref{eq:varphi}) and (\ref{eq:dt}) in Sect.~\ref{subsect:pfaff}.
The multi-valued function $U(t)$ is rewritten as  
\begin{align*}
U(t)&=
\prod_{k=1}^{L-1} \left(\frac{\tau_k}{\tau_0}\right)^{\zeta_k}  
\left(\frac{\tau_{k-1}}{\tau_0}-\frac{\tau_k}{\tau_0}\right)^{\eta_k}
\prod_{i=1}^{N} \left(1-x_i \frac{\tau_{L-1}}{\tau_0}\right)^{\theta_i}
\\
&={\tau_0}^{\zeta_0+1}
\prod_{k=1}^{L-1} {\tau_k}^{\zeta_k}  
(\tau_{k-1}-\tau_k)^{\eta_k}
\prod_{i=1}^{N} (\tau_0-x_i \tau_{L-1})^{\theta_i},
\end{align*}
where 
\[
\zeta_0 = -1 -\sum_{k=1}^{L-1}(\zeta_k+\eta_k)-\sum_{i=1}^N \theta_i .
\]

Now, we shall be concerned with  
the multi-valued function
\begin{align}
 \nonumber
 V(\tau)&=
 {\tau_0}^{-1} U(t)
 \\
 \label{eq:v}
&=
{\tau_0}^{\zeta_0}
\prod_{k=1}^{L-1} {\tau_k}^{\zeta_k}  
(\tau_{k-1}-\tau_k)^{\eta_k}
\prod_{i=1}^{N} (\tau_0-x_i \tau_{L-1})^{\theta_i}
\end{align}
and the rational $(L-1)$-forms
\begin{align}
\nonumber
\psi_0&=\tau_0 \varphi_0
=
\frac{\omega}{ \prod_{k=1}^{L-1}(\tau_{k-1}-\tau_k)},
\\
\label{eq:psi}
\psi_n^{(i)}&=\tau_0 \varphi_n^{(i)}
=
\frac{\omega}{ (x_i \tau_{L-1}-\tau_0)\prod_{\begin{subarray}{l} k=1 \\ k \neq n \end{subarray}}^{L-1}(\tau_{k-1}-\tau_k)}
\quad 
\left(
\begin{array}{c}
 1 \leq i \leq N \\
1 \leq n \leq L-1
\end{array}
\right)
\end{align}
in $\tau=(\tau_0,\ldots,\tau_{L-1})$.
Accordingly,
the integrals (\ref{eq:int}) can be expressed  
as
\[
y_0=\int_{\Delta} V(\tau) \psi_0,
\quad
y_n^{(i)}=\int_{\Delta} V(\tau) \psi_n^{(i)}
\]
since $U(t)\varphi_0=V(\tau)\psi_0$ and $U(t)\varphi_n^{(i)}=V(\tau)\psi_n^{(i)}$.
Without fear of repetition,
we summarize 
the correspondence  
\[
\zeta_k=
\alpha_k- \gamma_{k+1}, \quad 
\eta_k=\gamma_k-\alpha_k,
\quad
\theta_i=-\beta_i
\]
of constant parameters,  where
\[
\alpha_0= \sum_{i=1}^N \beta_i, \quad
\gamma_L=1.
\]

Next, let us introduce a `cyclic' transformation $\pi_i$
of the variables $\tau=(\tau_0,\ldots,\tau_{L-1})$
defined by
\[\pi_i: \tau_k \mapsto 
\left\{  \begin{array}{ll}
\tau_{k+1} & (0 \leq k \leq L-2)
\\ 
\tau_0/x_i & (k=L-1)
\end{array}
\right.
\]
for each $i=1,2, \ldots,N$.
We then have
\begin{align*}
{\pi_i}^n 
\left(\prod_{k=1}^{L-1} (\tau_{k-1}-\tau_k) \right)
&=\prod_{k=1}^{L-1} (\tau_{k+n-1}-\tau_{k+n}) 
\\
&=  (\tau_{n}-\tau_{n+1}) \cdots  (\tau_{L-2}-\tau_{L-1}) 
\underbrace{ (\tau_{L-1}-\tau_{L}) \cdots  (\tau_{L+n-2}-\tau_{L+n-1})}_{n}
\\
 &={x_i}^{-n}
 (x_i \tau_{L-1}-\tau_0)
 \prod^{L-1}_{ \begin{subarray}{c}
 k=1 \\
 k \neq n
 \end{subarray}}(\tau_{k-1}-\tau_k),
\end{align*}
where we regard  $\tau_{k}$ for $k \geq L$
as $\tau_k=\tau_{k-L}/x_i$
tentatively.
Combining this with $\pi_i(\omega)=(-1)^{L-1}\omega/x_i$,
we have
\begin{equation}
\label{eq:pipsi0}
{\pi_i}^n(\psi_0)=(-1)^{n(L-1)} \psi^{(i)}_{n}. 
\end{equation}
Note that 
the suffix $n$ of $\psi^{(i)}_{n}$ can be extended to be any $n \in {\mathbb Z}$ 
by the conditions
\[
\psi^{(i)}_{n+L}= \frac{ \psi^{(i)}_{n}  }{x_i} \quad \text{and} \quad
\psi^{(i)}_{0}=\psi_0.
\]
Hence we arrive at a `cyclic' expression 
\[
y_0=\int_{\Delta} V(\tau) \psi_0,
\quad
y_n^{(i)}=(-1)^{n(L-1)} \int_{\Delta} V(\tau) {\pi_i}^n(\psi_0)
\]
of the hypergeometric integrals (\ref{eq:int}).
Concerning the domain of integration, the $(L-1)$-simplex $\Delta_0$  can be written as
\[ \Delta_0=\{ 0 \leq \tau_{L-1} \leq \cdots \leq \tau_1 \leq \tau_0 \} \subset {\mathbb R}^L  \]
for instance; cf. (\ref{eq:Delta_0}).

\section{Case $N=1$: Thomae's ${}_L F_{L-1}$}
\label{sect:N=1}

This section is devoted to the case where $N=1$, i.e. the hypergeometric function
$F_{L,N}$ reduces to Thomae's ${}_LF_{L-1}$.
The linear Pfaffian system
under consideration is of the form
\begin{align*}
\tag{${\cal P}={\cal P}_{L,1}$}
\frac{{\rm d} \boldsymbol{y}}{{\rm d}x} 
&= \left\{
\frac{1}{x}
\begin{bmatrix}
0&
\\
a_1& b_1 
\\
a_2&a_2 & b_2
\\
\vdots&\vdots &\ddots &\ddots
\\
a_{L-1}&a_{L-1}&\cdots&a_{L-1}& b_{L-1}
\end{bmatrix}
+
\frac{1}{1-x}
\begin{bmatrix}
\beta&\beta& \cdots & \beta
\\
a_1& a_1 & \cdots &a_1
\\
a_2&a_2 & \cdots &a_2
\\
\vdots&\vdots &\ddots &\vdots
\\
a_{L-1}&a_{L-1}&\cdots&a_{L-1}
\end{bmatrix}
\right\}
\boldsymbol{y},
\\
\boldsymbol{y}&=
{}^{\rm T}\left(
y_0,y_1,\ldots,y_{L-1}\right),
\end{align*}
where $a_n=\alpha_n-\gamma_n$
and
$b_n(=b_{1,n})=-\gamma_n $ 
$(n=1,2,\ldots,L-1)$,
and its characteristic exponents at $x=0$ read $(0,b_1, \ldots,b_{L-1})$.

Our aim here is to write down the fundamental system of solutions $Y(x)$ of $\cal P$ having such a local power series expansion as
\begin{align}
 \nonumber
Y(x)&=
\Phi(x){\rm \, diag \, }
(1,x^{b_1}, \ldots,x^{b_{L-1}}),
\\
\label{eq:reqN=1}
\Phi(x)&=
\begin{bmatrix}
* &  &\\  
\vdots&\ddots&  \\
*&\cdots&*
\end{bmatrix} 
+O(x)=
 \left(  \boldsymbol{f}_0, \boldsymbol{f}_1,
\ldots,  \boldsymbol{f}_{L-1} \right),
\quad
\boldsymbol{f}_n
\in 
 {\mathbb C}^{L}[\![x]\!],
 \\  
 \nonumber
 \Phi(0)&\text{: invertible}
\end{align}
around the origin $x=0$.
For instance, 
we can take $ \boldsymbol{f}_0=\boldsymbol{y}(x;\Delta_0)$
for
$\Delta_0=\{0 \leq t_{L-1} \leq \cdots \leq t_2 \leq t_1 \leq 1 \} \subset {\mathbb R}^{L-1}$.

First we observe that 
the $n$th column $ \boldsymbol{f}_n$ $(n=1,2, \ldots,L-1)$ satisfies
the equation
\begin{equation}
\frac{{\rm d} \boldsymbol{f}_n}{{\rm d}x} 
= K_n \boldsymbol{f}_n,  
\label{eq:fn}
\end{equation}
where
\[
K_n=
\frac{1}{x}
\begin{bmatrix}
-b_n&
\\
a_1& b_1 -b_n
\\
a_2&a_2 & b_2 -b_n
\\
\vdots&\vdots &\ddots &\ddots
\\
a_{L-1}&a_{L-1}&\cdots&a_{L-1}& b_{L-1} -b_n
\end{bmatrix}
+
\frac{1}{1-x}
\begin{bmatrix}
\beta&\beta& \cdots & \beta
\\
a_1& a_1 & \cdots &a_1
\\
a_2&a_2 & \cdots &a_2
\\
\vdots&\vdots &\ddots &\vdots
\\
a_{L-1}&a_{L-1}&\cdots&a_{L-1}
\end{bmatrix}.
\]
Next we introduce a rotational matrix 
\[
\Lambda= 
\begin{bmatrix} 
0 & 1 && 
\\
 & 0 & \ddots &
\\
&&\ddots&1
\\
x^{-1} &&&0
\end{bmatrix}
\]
and consider
the equation satisfied by 
$\boldsymbol{g}= \Lambda^n \boldsymbol{f}_n$, namely,
\[
\frac{{\rm d} \boldsymbol{g}}{{\rm d}x} 
= \left(
\Lambda^n K_n \Lambda^{-n} + \frac{{\rm d} \Lambda^n}{{\rm d} x} \Lambda^{-n} \right)
\boldsymbol{g}.
\]
It is easy to verify the  following.

\begin{lemma} \label{lemma:nesw}
{\rm (1)} For a square matrix M of size $L$, it holds that
\[
\Lambda^n M \Lambda^{-n}
= \Lambda^n 
\left[
\begin{array}{c|c}
{\cal N} & {\cal E} \\
\hline
{\cal W} & {\cal S}
\end{array}
\right]
 \Lambda^{-n} 
 =
 \left[
\begin{array}{c|c}
{\cal S} & x {\cal W} \\
\hline
x^{-1} {\cal E} &{\cal N}
\end{array}
\right]
\]
with ${\cal N}$ being a square matrix of size $n$.
\\
{\rm(2)}
It holds that
\[
\frac{{\rm d} \Lambda^n}{{\rm d} x} \Lambda^{-n} 
=-\frac{1}{x}
 \left[
\begin{array}{c|c}
0 & 0 \\
\hline
0 & I_n
\end{array}
\right],
\]
where $I_n$ denotes the identity matrix of size $n$.
\end{lemma}

If we write the coefficient matrix $K_n$ of (\ref{eq:fn}) as
\[
K_n= \frac{1}{x}\left[
\begin{array}{c|c}
{\cal N} & {\cal E} \\
\hline
{\cal W} & {\cal S}
\end{array}
\right]
+\frac{1}{1-x}
\left[
\begin{array}{c|c}
{\cal N'} & {\cal E'} \\
\hline
{\cal W'} & {\cal S'}
\end{array}
\right]
\]
then,
by virtue of Lemma~\ref{lemma:nesw} and
${\cal E}=0$ and ${\cal W}={\cal W'}$, we have 
\begin{align*}
\widetilde{K}_n
&=
\Lambda^n K_n \Lambda^{-n} + \frac{{\rm d} \Lambda^n}{{\rm d} x} \Lambda^{-n}
\\
&=\frac{1}{x}
\left[
\begin{array}{c|c}
{\cal S}& x {\cal W} \\
\hline
x^{-1} {\cal E} & {\cal N}
\end{array}
\right]+
 \frac{1}{1-x}
\left[
\begin{array}{c|c}
{\cal S'}& x  {\cal W'} \\
\hline
x^{-1} {\cal E'} & {\cal N'}
\end{array}
\right]
-\frac{1}{x}
 \left[
\begin{array}{c|c}
0 & 0 \\
\hline
0 & I_n
\end{array}
\right]
\\
&=
\frac{1}{x}
\left[
\begin{array}{c|c}
{\cal S}& 0\\
\hline
{\cal E'} & {\cal N} -I_n
\end{array}
\right]+
 \frac{1}{1-x}
\left[
\begin{array}{c|c}
{\cal S'}& {\cal W'} \\
\hline
{\cal E'} & {\cal N'}
\end{array}
\right] .
\end{align*}
Hence, 
we obtain
\begin{align*}
\widetilde{K}_n
=&
\frac{1}{x}
\left[
\begin{array}{ccccc|cccc}
0 &&&&
\\
a_{n+1} & b_{n+1}-b_n &&&
\\
a_{n+2} &a_{n+2} &b_{n+2}-b_n &&
\\
\vdots & \vdots & \ddots& \ddots &
\\
a_{L-1} & a_{L-1} & \cdots & a_{L-1} & b_{L-1}-b_n &
\\
\hline
\beta &\beta &\cdots&\beta &\beta & -b_n-1
\\
a_1 &a_1 & \cdots&a_1 &a_1 &a_1 &b_1 -b_n-1
\\
\vdots & \vdots & & \vdots & \vdots &\vdots & \ddots & \ddots
\\
a_{n-1} &a_{n-1} & \cdots&a_{n-1} &a_{n-1} &a_{n-1} & \cdots &a_{n-1}& b_{n-1} -b_n-1
\end{array}
\right]
\\
&+ \frac{1}{1-x}
\left[
\begin{array}{cccc}
a_n & a_n & \cdots & a_n
\\
\vdots & \vdots & \ddots & \vdots
\\
a_{L-1}  & a_{L-1} &\cdots & a_{L-1}
\\
\hline
\beta  &\beta &\cdots & \beta 
\\
a_1 & a_1 &\cdots & a_1
\\
\vdots & \vdots & \ddots & \vdots
\\
a_{n-1} & a_{n-1} &\cdots & a_{n-1}
\end{array}
\right].
\end{align*}
 For convenience,
 we shall extend the suffixes $k$ of constant parameters to be any integers
by $a_{k+L}=a_k$, $b_{k+L}=b_k-1$, $a_0=\beta$ and $b_0=0$.
We therefore find that 
$\widetilde{K}_n$ is of the same form as the coefficient matrix
of the original ${\cal P}$
except the replacement
\[
a_k \mapsto a_{k+n}, \quad 
b_k \mapsto b_{k+n}-b_n, \quad
\beta \mapsto a_n
\]
of constant parameters.
In terms of the parameters
$(\alpha,\beta,\gamma)$, 
this replacement amounts to the transformation
\[
T_n: 
\alpha_k \mapsto \alpha_{k+n}-\gamma_n, \quad
(\alpha_0=) \  \beta \mapsto  \alpha_n-\gamma_n,
\quad
\gamma_k  \mapsto \gamma_{k+n}-\gamma_n,
\]
where 
$\alpha_{k+L}=\alpha_k+1$, $\gamma_{k+L}=\gamma_k+1$, 
 $\alpha_0=\beta$ and   $\gamma_0=0$.

Applying $T_n$ to the multi-valued function
\[U_0=U(t)=(-1)^{-\beta}\prod_{k=0}^{L-1} {t_k}^{\alpha_{k}-\gamma_{k+1}} 
 (t_{k-1}-t_k)^{\gamma_k-\alpha_k}
\]
shows that 
\[
U_n=T_n(U)=
 (-1)^{\gamma_n-\alpha_n}
\prod_{k=0}^{L-1} {t_k}^{\alpha_{k+n}-\gamma_{k+n+1}} 
 (t_{k-1}-t_k)^{\gamma_{k+n}-\alpha_{k+n}}.
\]
We shall write
\[
\varphi_n=\varphi^{(1)}_n=
\frac{\underline{{\rm d}t}}{\prod_{\begin{subarray}{l}
k=0 \\ k \neq n
\end{subarray}}^{L-1}(t_{k-1}-t_k)}
\]
with $t_0=1$ and $t_{-1}=x t_{L-1}$, and extendedly use the symbol 
 $\varphi_n$ for any integer $n$ by the quasi-periodicity 
 $\varphi_{n+L}=\varphi_{n}/x$. 
Now, the result can be  stated as follows.

\begin{thm}  \label{thm:rep1}
Let
\begin{align*}
\Phi(x)&=
\int_{\Delta_0}
\begin{bmatrix}
\varphi_0 & x \varphi_{L-1} & x \varphi_{L-2}  & \cdots & x  \varphi_{1} 
\\
\varphi_1  & \varphi_0 &x \varphi_{L-1} & \cdots & x  \varphi_{2} 
\\ 
\varphi_2  & \varphi_1 & \varphi_0  & \ddots & \vdots
\\
\vdots &\vdots &\ddots & \ddots & x\varphi_{L-1} 
\\
\varphi_{L-1} &\varphi_{L-2} & \cdots & \varphi_1& \varphi_0 
\end{bmatrix}
{\rm \, diag \,}
(U_0, U_1, \ldots,U_{L-1})
\\
&=\int_{\Delta_0} (U_n \varphi_{m-n})_{0 \leq m,n \leq L-1}.
\end{align*}
Then,
$Y(x)=\Phi(x){\rm \, diag \, }
(1,x^{-\gamma_1}, \ldots,x^{-\gamma_{L-1}})$
is a fundamental system of solutions of ${\cal P}(={\cal P}_{L,1})$
that fulfills the requirement 
{\rm(\ref{eq:reqN=1})}. 
\end{thm}

\noindent
We used only integrals over a single domain
$\Delta_0=\{0 \leq t_{L-1} \leq \cdots \leq t_2 \leq t_1 \leq 1 \}$
in the above theorem.

Furthermore, let us present an alternative expression of the same solution.
Introduce new variables $\tau=(\tau_0,\ldots,\tau_{L-1})$ and set
$t_n=\tau_n/\tau_0$.
Recall Sect.~\ref{sect:int}.
We thus have
$\Delta_0=\{0 \leq \tau_{L-1} \leq \cdots \leq \tau_{1} \leq \tau_0  \}$.
Define 
\[
\Delta_n= \pi^n(\Delta_0)=\{0 \leq \tau_{L-1+n} \leq \cdots \leq \tau_{1+n} \leq \tau_n  \}
\quad
\text{for $n \in{\mathbb Z}$}
\]
by using 
the `cyclic' transformation
$\pi: 
\tau_k \mapsto \tau_{k+1}$
with $\tau_{k+L}= {\tau_k}/{x}$.

\begin{thm}    \label{thm:rep2}
The $L \times L$ matrix function
\[
Y(x)
= \left(\int_{\Delta_n} U(t) \varphi_m \right)_{0 \leq m,n \leq  L-1}
\]
is a fundamental system of solutions of ${\cal P}$
that fulfills the requirement 
{\rm(\ref{eq:reqN=1})}. 
\end{thm}

\pf
Theorem~\ref{thm:rep1} implies that
$Y(x)$ can be expressed as
\begin{equation}\label{eq:rep1}
\int_{\Delta_0} (U_n \varphi_{m-n}x^{-\gamma_n})_{0 \leq m,n \leq L-1}
\times C
,
\end{equation}
where $\gamma_0=0$ and $C$ is a constant diagonal matrix.
We set
\begin{align*}
V(\tau)&={\tau_0}^{-1}U(t)
=(-1)^{-\beta} \prod_{k=0}^{L-1} 
{\tau_k}^{\alpha_{k}-\gamma_{k+1}} 
 (\tau_{k-1}-\tau_k)^{\gamma_k-\alpha_k},
 \\
\psi_n &= \tau_0 \varphi_n
=\frac{\omega}{\prod_{\begin{subarray}{l} k=0 \\ k \neq n \end{subarray}}^{L-1}(\tau_{k-1}-\tau_k)}
\end{align*}
and 
$\omega=\sum_{n=0}^{L-1} (-1)^n 
\tau_n
 {\rm d}\tau_0 \wedge \cdots \wedge  \widehat{{\rm d}\tau_n}
 \wedge \cdots\wedge {\rm d} \tau_{L-1}$
as well as in Sect.~\ref{sect:int}.
It holds that
$\pi^{-n}(V)
=x^{-\gamma_n}  (-1)^{\alpha_n-\gamma_n-\beta} {\tau_0}^{-1} U_n$
and
$\pi^{-n}(\psi_m)=(-1)^{n(L-1)} \psi_{m-n}$, 
where $\psi_{k+L}=\psi_{k}/x$.
Therefore,
\[
\pi^{-n}(U\varphi_m) =
\pi^{-n}(V \psi_m )= 
 x^{-\gamma_n}
 (-1)^{n(L-1)+\alpha_n-\gamma_n-\beta} U_n \varphi_{m-n},
 \]
 which completes the proof in view of (\ref{eq:rep1}).
\qed

\begin{example}[Gau\ss's ${}_2F_1$]
\rm

Let us restrict ourselves to the case where $(L,N)=(2,1)$, i.e.
the hypergeometric equation ${\cal P}_{L,N}$ thus becomes 
\begin{equation}
\label{eq:gauss}
\frac{{\rm d} \boldsymbol{y}}{{\rm d}x}
=\left(
\frac{1}{x}
\begin{bmatrix}
0 &0 \\
\alpha-\gamma & -\gamma
\end{bmatrix}
+ \frac{1}{1-x} 
\begin{bmatrix}
\beta &\beta \\
\alpha-\gamma & \alpha-\gamma
\end{bmatrix}
\right) \boldsymbol{y}.
\end{equation}
The multi-valued function $U=U(t)$ in a single variable $t=t_1$ and the rational $1$-forms $\varphi_n$ 
read
\[
U(t)=
 t^{\alpha-1} 
(1-t)^{\gamma-\alpha}
 (1-x t)^{-\beta},
 \quad
 \varphi_0= \frac{{\rm d} t}{1-t},
 \quad
  \varphi_1= \frac{{\rm d} t}{x t-1}.
\]
The domains
$\Delta_0=\{
0 \leq \tau_1 \leq \tau_0
\}$ and $\Delta_1=\{
0 \leq \tau_2 \leq \tau_1
\}$
of integration
are translated as
\[
\Delta_0=\{
0 \leq t \leq 1
\} \quad \text{and} \quad
\Delta_1=\{
1/x \leq t \leq \infty
\}
\]
by the correspondence $t_n=\tau_n/\tau_0$
and $\tau_{n+2}=\tau_n/x$.
Hence, it follows from Theorem~\ref{thm:rep2} 
that  the hypergeometric integrals
\begin{align*}
&
Y_{00}(x)=\int_0^1 {t}^{\alpha-1} (1-t)^{\gamma-\alpha-1}
(1-x t)^{-\beta} {\rm d} t,
\quad
Y_{01}(x)=\int_{1/x}^\infty {t}^{\alpha-1} (1-t)^{\gamma-\alpha-1}
(1-x t)^{-\beta} {\rm d} t,
\\
&Y_{10}(x)=-\int_0^1 {t}^{\alpha-1} (1-t)^{\gamma-\alpha}
(1-x t)^{-\beta-1} {\rm d} t,
\quad
Y_{11}(x)=-\int_{1/x}^\infty {t}^{\alpha-1} (1-t)^{\gamma-\alpha}
(1-x t)^{-\beta-1} {\rm d} t
\end{align*}
provide a fundamental system of solutions of (\ref{eq:gauss})
with the local behavior 
\[
Y(x) = \begin{bmatrix}
Y_{00}(x) &Y_{01}(x)\\
Y_{10}(x) &Y_{11}(x)
\end{bmatrix}
= \left(\begin{bmatrix}
* & 0 \\ * & *
\end{bmatrix}
+ O(x) \right) {\rm diag \ } (1, x^{-\gamma})
\]
near $x=0$.
\end{example}

\section{General $(L,N)$ case}
\label{sect:general}

Let us consider the domains 
\[\Delta_n^{(i)}={\pi_i}^n(\Delta_0)
 \quad
\left(
\begin{array}{c}
 1 \leq i \leq N \\
1 \leq n \leq L-1
\end{array}
\right)
\]
of integration,
defined by
applying the `cyclic' transformation
\[\pi_i: \tau_k \mapsto 
\left\{  \begin{array}{ll}
\tau_{k+1} & (0 \leq k \leq L-2)
\\ 
{x_i}^{-1} \tau_0 & (k=L-1)
\end{array}
\right. 
\]
to 
$\Delta_0=\{0 \leq \tau_{L-1} \leq \cdots \leq \tau_{1} \leq \tau_0  \} \subset {\mathbb R}^L$.

\begin{thm}
\label{thm:general}

The $(N(L-1)+1)\times(N(L-1)+1)$ 
matrix function 
\begin{align*}
Y(x)
&= \left(  \int_{\Delta_n^{(j)}} U\varphi_m^{(i)}
\right)_{ \begin{subarray}{l} 1 \leq i ,j \leq N \\ 0 \leq  m,n \leq L-1 \end{subarray} }
\\
& \text{\footnotesize \hspace{14mm} $0$ 
 \hspace{32mm} $1$ 
  \hspace{58mm} $N$}
\\
&=
\left[
\begin{array}{c|ccc|c|ccc}
\int_{\Delta_0}U\varphi_0 &\int_{\Delta_1^{(1)}}U\varphi_0 & \cdots &
\int_{\Delta_{L-1}^{(1)}}U\varphi_0 & \cdots &
\int_{\Delta_1^{(N)}}U\varphi_0 & \cdots &
\int_{\Delta_{L-1}^{(N)}}U\varphi_0 
\\
\hline
\int_{\Delta_0}U\varphi_1^{(1)} &\int_{\Delta_1^{(1)}}U\varphi_1^{(1)} & \cdots &
\int_{\Delta_{L-1}^{(1)}}U\varphi_1^{(1)} &  &
\int_{\Delta_1^{(N)}}U\varphi_1^{(1)} & \cdots &
\int_{\Delta_{L-1}^{(N)}}U\varphi_1^{(1)}
\\
\vdots & \vdots& \ddots &\vdots& \cdots & \vdots& \ddots &\vdots
\\
\int_{\Delta_0}U\varphi_{L-1}^{(1)} &\int_{\Delta_1^{(1)}}U\varphi_{L-1}^{(1)} & \cdots &
\int_{\Delta_{L-1}^{(1)}}U\varphi_{L-1}^{(1)} &  &
\int_{\Delta_1^{(N)}}U\varphi_{L-1}^{(1)} & \cdots &
\int_{\Delta_{L-1}^{(N)}}U\varphi_{L-1}^{(1)}
\\
\hline
\vdots & & \vdots && \ddots &&\vdots& \\
\hline
\int_{\Delta_0}U\varphi_1^{(N)} &\int_{\Delta_1^{(1)}}U\varphi_1^{(N)} & \cdots &
\int_{\Delta_{L-1}^{(1)}}U\varphi_1^{(N)} &  &
\int_{\Delta_1^{(N)}}U\varphi_1^{(N)} & \cdots &
\int_{\Delta_{L-1}^{(N)}}U\varphi_1^{(N)}
\\
\vdots & \vdots& \ddots &\vdots& \cdots & \vdots& \ddots &\vdots
\\
\int_{\Delta_0}U\varphi_{L-1}^{(N)} &\int_{\Delta_1^{(1)}}U\varphi_{L-1}^{(N)} & \cdots &
\int_{\Delta_{L-1}^{(1)}}U\varphi_{L-1}^{(N)} &  &
\int_{\Delta_1^{(N)}}U\varphi_{L-1}^{(N)} & \cdots &
\int_{\Delta_{L-1}^{(N)}}U\varphi_{L-1}^{(N)}
\end{array}
\right]
\end{align*}
is a fundamental system of solutions of ${\cal P}_{L,N}$ 
having the local behavior stated in
 Theorem~\ref{thm:main}.
\end{thm}

\pf
Fix $j \in \{ 1, \ldots, N \}$ and $n \in \{1, \ldots,L-1\}$.
We shall
examine the local behavior of the 
column vector 
 \begin{equation}  \label{eq:vecy}
\boldsymbol{y}(x; \Delta^{(j)}_n) =\left(\int_{\Delta_n^{(j)}} U\varphi_m^{(i)} \right)_{ \begin{subarray}{l} 1 \leq i \leq N  \\ 0 \leq  m \leq L-1 \end{subarray} }
 \end{equation}
 of
hypergeometric integrals,
which belongs to the $j$th block of the matrix function $Y(x)$.
To this end, we first rewrite its element as
an integral over the simplex $\Delta_0$, namely
\[
\int_{\Delta_n^{(j)}} U\varphi_m^{(i)} 
=\int_{\Delta_n^{(j)}} V\psi_m^{(i)}
 =\int_{\Delta_0} {\pi_j}^{-n} \left(V\psi_m^{(i)} \right).
\]
Recall (\ref{eq:v}) and (\ref{eq:psi}) for notations.
Here, we mention that
$\boldsymbol{y}(x; \Delta^{(j)}_n)$ certainly solves ${\cal P}_{L,N}$
since 
$\Delta_n^{(j)}$ is a chamber framed by the hyperplanes
which are the singular loci
of the multi-valued function $U=U(t)$;
cf. \cite{tsu12}.

Applying ${\pi_j}^{-n}$ to
$V=V(\tau)$ yields
\begin{align*}
{\pi_j}^{-n}(V)&={\tau_{-n}}^{\zeta_0}
\prod_{k=1}^{L-1} {\tau_{k-n} }^{\zeta_k}  
(\tau_{k-n-1}-\tau_{k-n})^{\eta_k}
\prod_{i=1}^{N} (\tau_{-n}-x_i \tau_{L-n-1})^{\theta_i}
\\
&=
{x_j}^{ \zeta_0+\sum_{k=1}^{n-1} (\zeta_k + \eta_k) + \theta_j}
 {\tau_{L-n}}^{\zeta_0}
 \prod_{k=1}^{n-1}
{\tau_{L+k-n} }^{\zeta_k}  
(\tau_{L+k-n-1}-\tau_{L+k-n})^{\eta_k}
\\
&\quad
\times
 \prod_{k=n}^{L-1} {\tau_{k-n} }^{\zeta_k}  
(\tau_{k-n-1}-\tau_{k-n})^{\eta_k}
\\
& \quad
\times (\tau_{L-n}-\tau_{L-n-1})^{\theta_j}
 \prod_{i \neq j} (x_j\tau_{L-n}-x_i\tau_{L-n-1})^{\theta_i}
\\
&=:{x_j}^{b_{j,n}} g(x, \tau)
\end{align*}
by the use of  the quasi-periodicity
$\tau_{k+L}=\tau_k/x_j$ and
$b_{j,n}=\sum_{i \neq j} \beta_i -\gamma_n =
\zeta_0
+\sum_{k=1}^{n-1} (\zeta_k + \eta_k) + \theta_j$.
Observe that $g(x,\tau)$ is holomorphic on the interior
${\rm Int}(\Delta_0)$
of $\Delta_0$
provided 
\begin{equation}
\label{eq:open}
|x_i|>|x_j| \quad \text{for} \quad i \neq j.
\end{equation}
It readily follows from (\ref{eq:pipsi0}) that
\[
{\pi_j}^{-n}(\psi^{(j)}_m)
= \left\{
\begin{array}{ll}
(-1)^{n(L-1)} x_j \psi^{(j)}_{L+m-n}  & (m < n )
\\
(-1)^{n(L-1)} \psi_0 & (m=n)
\\
(-1)^{n(L-1)} \psi^{(j)}_{m-n}  & (m > n )
\end{array}
\right.
\]
and thus
${\pi_j}^{-n}(\psi^{(j)}_m)$ is holomorphic on 
${\rm Int}(\Delta_0)$
as long as $|x_j|<1$.
If we remember that any solution of the linear Pfaffian system
${\cal P}_{L,N}$
is holomorphic outside its singular locus $\Xi$ (see (\ref{eq:Xi}))
and notice that (\ref{eq:open}) is an open condition and  thereby removable via the identity theorem,
then we can conclude
 that
\[
\int_{\Delta_0} {\pi_j}^{-n} \left(V\psi_m^{(j)} \right)
= \left\{
\begin{array}{ll}
{x_j}^{b_{j,n}+1} \times (\text{holomorphic on $D^{(j)}$})   & (m < n )
\\
{x_j}^{b_{j,n}}  \times (\text{holomorphic on $D^{(j)}$})     & (m \geq n)
\end{array}
\right.
\]
where
$D^{(j)} =\{ |x_j|<1  \} \cap \bigcap_{i \neq j} \{ x_i \neq x_j\} \subset {\mathbb C}^N$.

Next we deal with the case where $i \neq j$. 
Applying ${\pi_j}^{-n}$ to the  denominator and numerator of 
\[
\psi_m^{(i)}=\frac{\omega}{ (x_i \tau_{L-1}-\tau_0)\prod_{\begin{subarray}{l} k=1 \\ k \neq m \end{subarray}}^{L-1}(\tau_{k-1}-\tau_k)}
\]
yields
\begin{align*}
\lefteqn{
{\pi_j}^{-n}
\left(
(x_i \tau_{L-1}-\tau_0)
\prod_{\begin{subarray}{l}
 k=1 \\
 k \neq m
 \end{subarray}}^{L-1} (\tau_{k-1}-\tau_k)
  \right)
  }
  \\
  &=
(x_i \tau_{L-n-1}-\tau_{-n})
\prod_{\begin{subarray}{l}
 k=1 \\
 k \neq m
 \end{subarray}}^{L-1} (\tau_{k-n-1}-\tau_{k-n})
 \\
 &=(x_i \tau_{L-n-1}- x_j\tau_{L-n})
\underbrace{ (\tau_{-n}-\tau_{1-n}) \cdots (\tau_{-2}-\tau_{-1}) }_{n-1}
 (\tau_{-1}-\tau_{0}) \cdots  (\tau_{L-n-2}-\tau_{L-n-1})  
 \\
 &\quad
 \times \frac{1}{\tau_{m-n-1}-\tau_{m-n}}
 \\
 &= \left\{
 \begin{array}{ll}
{x_j}^{n-2} (x_i \tau_{L-n-1}- x_j\tau_{L-n})
(x_j \tau_{L-1}-\tau_{0}) 
\prod^{L-1}_{\begin{subarray}{l}
 k=1 \\
 k \neq L+m-n, L-n
 \end{subarray}  }
 (\tau_{k-1}-\tau_k)
& (m <n)
 \\
 {x_j}^{n-1} (x_i \tau_{L-n-1}- x_j\tau_{L-n}) 
 \prod^{L-1}_{\begin{subarray}{l}
 k=1 \\
 k \neq L-n
 \end{subarray}  }
 (\tau_{k-1}-\tau_k)
 &(m=n)
 \\
 {x_j}^{n-1}  (x_i \tau_{L-n-1}- x_j\tau_{L-n})
 (x_j \tau_{L-1}-\tau_{0}) 
\prod^{L-1}_{\begin{subarray}{l}
 k=1 \\
 k \neq m-n, L-n
 \end{subarray}  }
 (\tau_{k-1}-\tau_k)
 & (m >n)
 \end{array}
 \right.
\end{align*}
and
${\pi_j}^{-n}(\omega)=(-1)^{n(L-1)}{x_j}^n \omega$.
Therefore, 
taking  an integral of ${\pi_j}^{-n}(V \psi_m^{(i)})$ over $\Delta_0$
shows that
\[
\int_{\Delta_0} {\pi_j}^{-n} \left(V\psi_m^{(i)} \right)
= \left\{
\begin{array}{ll}
{x_j}^{b_{j,n}+2} \times (\text{holomorphic on $D^{(j)}$})   & (m < n )
\\
{x_j}^{b_{j,n}+1}  \times (\text{holomorphic on $D^{(j)}$})     & (m \geq n)
\end{array}
\right.
\]
for $i \neq j$
in the same manner as above.

We have verified that
(\ref{eq:vecy}) certainly possesses the characteristic behavior 
specified in Theorem~\ref{thm:main}. 
\qed

\section{From hypergeometric equation ${\cal P}_{L,N+1}$ to isomonodromic deformations}

\label{sect:mpd}

The subject of this section is a connection between 
the hypergeometric equation and isomonodromic deformations 
of a certain Fuchsian system,
from which the hypergeometric solution of the Painlev\'e 
equation naturally arises 
as a by-product; cf. \cite{tsu12}.

Consider the linear Pfaffian system
${\cal P}_{L,N+1}$
of rank $(N+1)(L-1)+1$.
Suppose $\beta_{N+1}=0$.
It is then obvious that $U=U(t)$ does not depend on 
$x_{N+1}$; see (\ref{eq:U}).
Accordingly, 
the $N(L-1)+1$ functions
\[
y_0=\int_\Delta U(t) \varphi_0, \quad y_{n}^{(i)}=\int_\Delta U(t) \varphi_{n}^{(i)} \quad
\left(
\begin{array}{c}
 1 \leq i \leq N \\
1 \leq n \leq L-1
\end{array}
\right)
\] 
do not  depend on  $x_{N+1}$ 
and thus constitute a solution of ${\cal P}_{L,N}$
if the domain $\Delta$ of integration is suitably chosen. 
We are now interested in how 
the other $L-1$ functions $y^{(N+1)}_n$ $(1 \leq n \leq L-1)$
depend on $x_{N+1}$.

We take the change of variables
\[
 x_i=\frac{1}{u_i}  \quad(1 \leq i \leq N), \quad  x_{N+1}=\frac{1}{z}
 \]
 and rewrite the constant parameters as
\begin{equation} \label{eq:param}
\alpha_n=e_n-e_0, \quad 
\beta_i=  - \theta_i, \quad  
\gamma_n= e_n-e_0-\kappa_n
\quad
\left(
\begin{array}{c}
 1 \leq i \leq N \\
1 \leq n \leq L-1
\end{array}
\right)
\end{equation}
for the sake of convenience.
Let
\[
f_0=1, \quad f_n= -y_n^{(N+1)}\prod_{j=1}^N {u_j}^{\theta_j} \quad (1\leq n \leq L-1).
\]
Then, we verify from ${\cal P}_{L,N+1}$ that
$\boldsymbol{f}={}^{\rm T}\left( f_0,f_1,\ldots,f_{L-1}\right)$ 
satisfies the Fuchsian system 
\begin{equation} \label{eq:A'}
\frac{{\rm d}\boldsymbol{f}}{{\rm d} z}
= A'  \boldsymbol{f}
= \sum_{i=0}^{N+1} \frac{A_i'}{z-u_i} \boldsymbol{f}
\quad 
(u_0=1, \ u_{N+1}=0)
\end{equation}
of ordinary differential equations with respect to $z=1/x_{N+1}$, 
whose coefficients are given by
\begin{align*}
A_0' &= 
 \begin{bmatrix}
 0 & 0 & \cdots & 0
 \\
\kappa_1 h& - \kappa_1 & \cdots  & - \kappa_1
\\
\vdots & \vdots & \ddots & \vdots
\\
 \kappa_{L-1} h & - \kappa_{L-1} & \cdots  & - \kappa_{L-1}
\end{bmatrix}, 
\quad
A_i' 
=
\begin{bmatrix}
 0 & 0 & \cdots & 0
 \\
b_1^{(i)}& \theta_i & \cdots  & 0 
\\
\vdots & \vdots & \ddots & \vdots
\\
b_{L-1}^{(i)} & 0& \cdots  & \theta_i
\end{bmatrix}
\quad (1 \leq i \leq N)
\\
A_{N+1}' 
&= 
 \begin{bmatrix}
   0 &   0   & 0  & 0 & \cdots  & 0
   \\
        & e_1-e_0 & \kappa_1 & \kappa_1 &  \cdots     &  \kappa_1
   \\     
          &        &e_2-e_0& \kappa_2 &  \cdots     &  \kappa_2
   \\
          &        &          & e_3 -e_0&  \ddots     &  \vdots
   \\
          &         &          &           & \ddots  & \kappa_{L-2}
    \\
          &         &          &            &              & e_{L-1}-e_0
\end{bmatrix}
\end{align*}
with
\[
h= y_0 \prod_{j=1}^N {u_j}^{\theta_j}, \quad 
b_n^{(i)} =y^{(i)}_n
\theta_i\prod_{j=1}^N {u_j}^{\theta_j} 
\quad 
\left(
\begin{array}{c}
 1 \leq i \leq N \\
1 \leq n \leq L-1
\end{array}
\right).
\]
It goes without saying that 
(\ref{eq:A'}) is reducible.
More properly,
its solution space turns out to be
a direct sum of ${\mathbb C} \boldsymbol{f}$ 
and the solution space of
${\cal P}_{L-1,1}$;
see Remark~\ref{remark:thomae} below.

Let $Y'$ denote the fundamental system of solutions of (\ref{eq:A'}).
Its gauge transformation
\[
Y= \left(
z^{e_0}
\prod_{i=1}^N  (z-u_i)^{-\theta_i}
{u_i}^{\theta_i/L} \right) Y'
\]
then satisfies the  Fuchsian system 
\begin{equation}  \label{eq:lax1}
\frac{{\rm d} Y}{{\rm d} z}
= A
Y
= \sum_{i=0}^{N+1} \frac{A_i}{z-u_i} Y
\end{equation}
whose coefficients are given by 
$A_0=A_0'$, $A_i=A_i'-\theta_i I_L$ $(1 \leq i \leq N)$
and $A_{N+1}=A_{N+1}'+e_0 I_L$.
We know {\it a priori} that the monodromy of $Y(z)$ is independent of 
variables $u_i$ $(1 \leq i \leq N)$;
see Remark~\ref{remark:rigid}.
Actually, we verify again from ${\cal P}_{L,N+1}$ that 
$Y$ satisfies the extended system
\begin{equation}   \label{eq:lax2}
\frac{\partial Y}{\partial u_i}
= B_i Y
\quad (1 \leq i \leq N)
\end{equation}
of linear differential equations with respect to
$u_i$,
whose coefficients are given by
\[
B_i
=
\frac{1}{u_i-z}
\begin{bmatrix}
- \theta_i & 0  & \cdots& 0
 \\
b_1^{(i)}& & &
\\
\vdots& &\mbox{\huge 0}&
\\
b_{L-1}^{(i)}&  &  & 
\end{bmatrix}
-
\frac{1}{u_i} 
\begin{bmatrix}
 -\frac{\theta_i}{L}& & & 
 \\
b_1^{(i)}& -\frac{\theta_i}{L}&&
\\
\vdots&&\ddots&
\\
b_{L-1}^{(i)}&&& -\frac{\theta_i}{L}
\end{bmatrix}.
\]
The coefficient $B_i$ is a rational function in $z$ and, thus, 
 (\ref{eq:lax2}) describes the isomonodromic family of  (\ref{eq:lax1})
 along the deformation parameters $u_i$.
 Namely,  (\ref{eq:lax2}) guarantees the existence of a fundamental system of solutions of 
 (\ref{eq:lax1})
 whose monodromy matrices do not depend on $u_i$;
 see \cite{sch12}.

The Riemann scheme
of (\ref{eq:lax1})
is written as follows:
\[
\begin{array}{|c|c|}
\hline
\text{Singularity} & \text{Characteristic exponents}  \\ \hline
z=0 & (e_0,e_1,\ldots,e_{L-1})  \\ \hline    
z=\infty & (\kappa_{0}-e_0,\kappa_{1}-e_1, \ldots,\kappa_{L-1}-e_{L-1} )  \\ 
\hline
z=u_0=1 & \left(-\sum_{n=1}^{L-1} \kappa_n ,0, \ldots,0 \right)
\\
\hline
z=u_i  \mbox{\ } (1 \leq i \leq N )& (-\theta_i, 0, \ldots,0)  \\
\hline
   \end{array}
\]
(However, the relation 
$\kappa_0=\sum_{i=1}^N \theta_i$
holds.)
Accordingly,
its spectral type reads
\[
\begin{array}{ll}
\underbrace{1,1, \ldots,1}_{L}
&
\text{at $z=0, \infty$}, 
\\
1, L-1
&
\text{at $z=u_i$ $(0 \leq i \leq N )$} .
\end{array}
\]

In the {\it general} case, 
a Fuchsian system with the above spectral type is non-rigid
(cf. Remark~\ref{remark:rigid});
in fact, it is equipped with $2N(L-1)$ accessory parameters.
The coefficients of
 such a Fuchsian system 
 \[ \tag{${\cal L}_{L,N}$}
 \frac{{\rm d} Y}{{\rm d} z}
= \sum_{i=0}^{N+1} \frac{{\cal A}_i}{z-u_i} Y
\]
 can be parametrized as
 \begin{align*}
{\cal A}_i 
&=  {}^{\rm T}\left(b_0^{(i)}, b_1^{(i)}, \ldots ,b_{L-1}^{(i)}\right)
\cdot
\left(c_0^{(i)}, c_1^{(i)}, \ldots ,c_{L-1}^{(i)}\right)\quad
(0 \leq i \leq N),
\\
{\cal A}_{N+1} 
&= 
 \begin{pmatrix}
   e_0 &   w_{0,1}      & \cdots  & w_{0,L-1}
   \\
        & e_1 &  \ddots       & \vdots
   \\     
          &        &\ddots& w_{L-2,L-1}
   \\
          &        &          & e_{L-1}    
\end{pmatrix}
\quad \text{with} \quad 
w_{m,n}=-\sum_{i=0}^N b_m^{(i)} c_n^{(i)},
\end{align*}
where
$c_0^{(i)}=1$,
${\rm tr \,} {\cal A}_i =
  \sum_{n=0}^{L-1} b_n^{(i)}c_n^{(i)}
 = -\theta_i$
 and 
$ \sum_{i=0}^{N} b_n^{(i)} c_n^{(i)}=- \kappa_n$.
We can and will normalize
 the characteristic exponents
by
\begin{equation}
\label{eq:norm1}
{\rm tr \,} {\cal A}_{N+1}=\sum_{n=0}^{L-1} e_n=\frac{L-1}{2}
\end{equation}
without loss of generality.
Assume 
the Fuchsian relation
\begin{equation}
\label{eq:norm2}
\sum_{n=0}^{L-1}\kappa_n
=\sum_{i=0}^N \theta_i
\end{equation}
holds.
As shown in \cite{tsu13},
the isomonodromic deformations of 
${\cal L}_{L,N}$
are 
governed by the Hamiltonian system
\[
\tag{${\cal H}_{L,N}$}
\frac{\partial q_n^{(i)}}{\partial x_j}=\frac{\partial H_j}{\partial p_n^{(i)}},
\quad 
\frac{\partial p_n^{(i)}}{\partial x_j}=-\frac{\partial H_j}{\partial q_n^{(i)}}
\quad 
\left(
\begin{array}{c}
 1 \leq i , j \leq N \\
1 \leq n \leq L-1
\end{array}
\right)
\]
of partial differential equations with respect to 
variables $x_i=1/u_i$,
whose Hamiltonian function $H_i$ is defined by 
\[
x_i H_i=
\sum_{n=0}^{L-1} e_n q_n^{(i)} p_n^{(i)} 
+\sum_{j=0}^N \sum_{0 \leq m <n \leq L-1}
 q_m^{(i)}   p_m^{(j)} q_n^{(j)} p_n^{(i)} 
+\sum^N_{
\begin{subarray}{c} 
j=0  \\
j \neq i
\end{subarray}
}
\frac{x_j}{x_i-x_j} 
 \sum_{m,n=0}^{L-1}
 q_m^{(i)}   p_m^{(j)} q_n^{(j)} p_n^{(i)}
\]
with
$x_0=q_n^{(0)}
= q_0^{(i)}=1$, 
$p_n^{(0)}
=
\kappa_{n} - \sum_{i=1}^{N} q_n^{(i)} p_n^{(i)}$
and 
$p_0^{(i)}
=\theta_i-\sum_{n=1}^{L-1} q_n^{(i)} p_n^{(i)}$.
Thus, $H_i$ becomes a polynomial in the 
$2N(L-1)$
unknowns ({\it canonical coordinates})
\begin{equation}
\label{eq:coord}
q_n^{(i)}=\frac{c_n^{(i)}}{c_n^{(0)}} \quad \text{and} \quad 
p_n^{(i)}=-b_n^{(i)}c_n^{(0)}
\quad 
\left(
\begin{array}{c}
 1 \leq i \leq N \\
1 \leq n \leq L-1
\end{array}
\right).
\end{equation}
The Hamiltonian system ${\cal H}_{L,N}$ contains constant parameters
\[
(e,\kappa,\theta)=(e_0,\ldots,e_{L-1},\kappa_0,\ldots,\kappa_{L-1},\theta_0,\ldots,\theta_N),
\]
but their number is essentially 
 $2L+N-1$
 on account of 
(\ref{eq:norm1}) and (\ref{eq:norm2}).
Note that 
${\cal H}_{L,N}$ is an extension of  the sixth Painlev\'e equation,
as it literally recovers the original if $(L,N)=(2,1)$.

Interestingly enough, the previous system (\ref{eq:lax1}) coincides with
a {\it particular} case of ${\cal L}_{L,N}$ such that
\begin{align*}
&\kappa_0=\sum_{i=1}^N \theta_i 
\quad \text{and}
\\
&b_0^{(0)}=0, \quad b_0^{(i)}=-\theta_i,
\quad
c_n^{(0)}=\frac{-1}{h}, \quad
c_n^{(i)}=0 \quad  
\left(
\begin{array}{c}
 1 \leq i \leq N \\
1 \leq n \leq L-1
\end{array}
\right).
\end{align*}
Hence, if we remember the algebraic relations (\ref{eq:coord})
between 
 the canonical coordinates of ${\cal H}_{L,N}$ and the coefficients of
its associated Fuchsian system ${\cal L}_{L,N}$,
then we derive directly from the above argument
an $N(L-1)$-parameter family of particular solutions of ${\cal H}_{L,N}$;
cf. \cite[Theorem 3.2]{tsu12}.

\begin{thm}
If $\kappa_0=\sum_{i=1}^N \theta_i$,
then the Hamiltonian system
${\cal H}_{L,N}$ 
has
a particular solution 
\[
q_n^{(i)}=0 \quad \text{and} \quad
p_n^{(i)}
=
\theta_i\frac{y_n^{(i)}}{y_0},
\]
where $\left\{y_0,y_n^{(i)}\right\}$ is an arbitrary solution of 
the linear Pfaffian system
${\cal P}_{L,N}$
with {\rm(\ref{eq:param})}.
\end{thm}

\begin{remark}\rm
\label{remark:thomae}

A fundamental system $Y'$ of solutions of (\ref{eq:A'}) can be taken of the form
\[
Y'= 
\begin{bmatrix}
 f_0& 0& \cdots & 0 
 \\
f_1&&&
\\
\vdots&&\mbox{\huge $W$}&
\\
f_{L-1}&&& 
\end{bmatrix}.
\]
The gauge transformation $Z= \left(z^{\kappa_1-e_1+e_0} \prod_{i=1}^N (z-u_i)^{-\theta_i} \right) W$  
satisfies the Fuchsian system
\begin{align*}
\frac{{\rm d} Z}{{\rm d}w}
&=\frac{1}{w}
 \begin{bmatrix}
  0   \\     
  \kappa_2  &\kappa_2-\kappa_1-e_2+e_1 \\
  \kappa_3  &  \kappa_3 & \kappa_3-\kappa_1-e_3 +e_1 \\
  \vdots &\vdots& \ddots & \ddots  & 
    \\
  \kappa_{L-1}& \kappa_{L-1}& \cdots & \kappa_{L-1}&  \kappa_{L-1}-\kappa_1- e_{L-1}+e_1
\end{bmatrix}
Z
\\
& \qquad 
+
\frac{1}{1-w}
\begin{bmatrix}
  \kappa_1 & \cdots  & \kappa_1
\\
 \vdots & \ddots & \vdots
\\
 \kappa_{L-1} & \cdots  & \kappa_{L-1}
\end{bmatrix}
Z
\end{align*}
with $w=1/z$,
which is exactly the hypergeometric equation ${\cal P}_{L-1,1}$
(for Thomae's ${}_{L-1}F_{L-2}$)
under the correspondence of constant parameters as
\begin{align*}
&\alpha_1=\kappa_1-e_1+e_2, \quad
\alpha_2=\kappa_1-e_1+e_3,
\quad \ldots ,   
\quad
\alpha_{L-2}=\kappa_1-e_1+e_{L-1}, \quad 
\beta= \kappa_1,
\\
&\gamma_1=\kappa_1-\kappa_2 -e_1+e_2, \quad
\gamma_2=\kappa_1-\kappa_3 -e_1+e_3, \quad \ldots ,   
\quad
\gamma_{L-2}=\kappa_1-\kappa_{L-1} -e_1+e_{L-1}.
\end{align*}
Cf. Sect.~\ref{subsect:pfaff}.
\end{remark}

\small
\paragraph{\it Acknowledgement.}
The initial idea of the content of Sect.~\ref{sect:mpd} is due to Yamada \cite{yam10}, which was reported to me a few days after I sent him the first draft of the article \cite{tsu12}.  
I would like to express my sincere gratitude to Yasuhiko Yamada
for his permission to use it.
I also thank the referees for reading carefully the original manuscript and giving valuable 
suggestions.

\end{document}